\documentclass{article}
\usepackage{amssymb}
\usepackage{graphics}
\input xy
\xyoption{all}

\begin{document}
\title{Structure of the curvature tensor\\ on symplectic spinors}
\author{Svatopluk Kr\'ysl \footnote{{\it E-mail address}: krysl@karlin.mff.cuni.cz}\\ {\it \small  Charles University of Prague, Sokolovsk\'a 83, Praha,  Czech Republic}.
\thanks{
 The author of this article was supported by the grant
GA\v{C}R 201/06/P223 of the Grant Agency of Czech Republic for young
researchers. The work is a part of the research project MSM 0021620839 financed by M\v{S}MT \v{C}R.}}

\maketitle \noindent
\centerline{\large\bf Abstract} We study symplectic manifolds $(M^{2l},\omega)$ equipped with a symplectic torsion-free 
affine (also called Fedosov) connection $\nabla$ and admitting a metaplectic structure. Let $\mathcal{S}$ be the so called
 symplectic spinor bundle and let $R^S$ be the curvature tensor field of the symplectic spinor covariant derivative $\nabla^S$ 
associated to the Fedosov connection $\nabla.$ It is known that the space of symplectic spinor valued exterior differential $2$-forms,  
$\Gamma(M,\bigwedge^2T^*M\otimes {\mathcal{S}}),$ decomposes into three invariant spaces with respect to the structure group, which is the metaplectic 
group $Mp(2l,\mathbb{R})$ in this case. For a symplectic spinor field $\phi \in \Gamma(M,\mathcal{S}),$ we compute explicitly the projections of 
$R^S\phi \in \Gamma(M,\bigwedge^2T^*M \otimes \mathcal{S})$ onto the three mentioned invariant spaces in terms of the symplectic Ricci and symplectic 
Weyl curvature tensor fields of the connection $\nabla.$ Using this decomposition, we derive a complex of first order differential operators provided the 
Weyl tensor of the Fedosov connection is trivial.

{\it Math. Subj. Class.:} 53C07, 53D05, 58J10.

{\it Keywords:} Fedosov manifolds, metaplectic structures, symplectic spinors, Kostant spinors, Segal-Shale-Weil representation, sympletic curvature tensor.

\section{Introduction}

In the paper, we shall study the action of the symplectic curvature tensor field on symplectic spinors over a symplectic manifold $(M^{2l},\omega)$ with a given metapletic structure and equipped with a symplectic torsion-free affine connection $\nabla$. Such connections are usually called  Fedosov connections.
It is well known that in the case of $l>1,$ the curvature tensor field of the connection $\nabla$ decomposes into two  parts, namely into the symplectic Weyl and the symplectic Ricci curvature tensor field. In the case $l=1,$ only the symplectic Ricci curvature tensor field appears. See Vaisman \cite{Vaisman} for details. 

Now, let us say a few words about the metaplectic structure.
In the symplectic case, there exists (in a parallel to the Riemannian case) a non-trivial two-fold covering of the symplectic group
$Sp(2l,\mathbb{R}),$ the so called metaplectic group. We shall denote it by $Mp(2l,\mathbb{R}).$ 
A metaplectic structure on a symplectic manifold $(M^{2l},\omega)$ is
a notion parallel to a spin structure on a Riemannian manifold. 
 For a symplectic manifold admitting a metaplectic structure, one can construct the so called {\it symplectic spinor bundle} $\mathcal{S},$  introduced by
  Bertram Kostant in 1974.  The symplectic spinor bundle $\mathcal{S}$ is the vector bundle associated to the metaplectic structure
  on $M$ via the so called  Segal-Shale-Weil representation
of the metaplectic group $Mp(2l,\mathbb{R})$. See Kostant \cite{Kostant2} for details.
 
 The Segal-Shale-Weil representation is an infinite dimensional unitary representation of
the metaplectic group $Mp(2l, \mathbb{R})$ on the space of all complex valued square
Lebesgue integrable functions ${\bf L^{2}}(\mathbb{R}^{l}).$ Because of the infinite dimension, the  Segal-Shale-Weil representation is not so is to handle.
It is known, see, e.g., Kashiwara, Vergne \cite{KV}, that  
the underlying Harish-Chandra module of
this representation is equivalent to the space
$\mathbb{C}[x^1,\ldots, x^l]$ of polynomials in $l$ variables, on which
the Lie algebra $\mathfrak{sp}(2l,\mathbb{C})$ acts 
via the so called Chevalley homomorphism,\footnote{The Chevalley homomorphism
realizes the complex symplectic Lie algebra as a Lie subalgebra of the
algebra of polynomial coefficients differential operators acting on
$\mathbb{C}[x^1,\ldots, x^l].$} see Britten, Hooper, Lemire
\cite{BHL}. Thus, the infinitesimal structure of the Segal-Shale-Weil
representation can be viewed as the complexified {\it symmetric}
 algebra $(\bigoplus_{i=0}^{\infty}\odot^i
\mathbb{R}^l)\otimes_{\mathbb{R}}{\mathbb{C}}$ of the Lagrangian
subspace $\mathbb{R}^l$ of the canonical symplectic vector space
$\mathbb{R}^{2l}\simeq \mathbb{R}^l\oplus \mathbb{R}^l.$
This shows that the situation is completely parallel to
the complex orthogonal case  and the spinor representation, which can be
realized as the {\it exterior} algebra of a maximal isotropic
subspace.  An interested reader is
referred to Weil \cite{Weil}, Kashiwara, Vergne \cite{KV} and also
to Britten, Hooper, Lemire \cite{BHL} for more details.  For some technical reasons, we shall be using the so called minimal globalization of
the Harish-Chandra $(\mathfrak{g},K)$-module of the  Segal-Shale-Weil representation, which we will call {\it metaplectic representation} and  denote it by ${\bf S}$ (the elements of $\bf S$ will be called {\it symplectic spinors}). This representation, as well as the Segal-Shale-Weil one, decomposes into two irreducible subrepresentations ${\bf S}_{+}$ and ${\bf S}_-.$  

For any symplectic connection $\nabla$ on a symplectic manifold $(M,\omega)$ admitting a metaplectic structure, we can form the associated covariant derivative $\nabla^S$ acting on the sections of the symplectic spinor bundle $\mathcal{S}$. 
The curvature tensor field $R^S: \Gamma(M,\mathcal{S}) \to \Gamma(M, \bigwedge^2 TM^* \otimes \mathcal{S})$ of the associated covariant derivative $\nabla^S$ is defined by the classical formula. The tensor field $R^S$ decomposes also into two parts, one of which is dependent only on the symplectic Ricci and the second one on the symplectic Weyl tensor. 
It is known (cf. Kr\'{y}sl \cite{KryslSVF}) that the space $\bigwedge^2 \mathbb{R}^{2l} \otimes {\bf S}_{\pm}$ decomposes into three irreducible summands wr. to the natural action of $Mp(2l,\mathbb{R})$ on this space.  We shall briefly describe this result in the paper.    Let us denote the mentioned three summands of the decomposition of $\bigwedge^2\mathbb{R}^{2l}\otimes {\bf S}_{\pm}$ by ${\bf E}^{20}_{\pm},$ ${\bf E}^{21}_{\pm}$ and ${\bf E}^{22}_{\pm}$ and the corresponding vector bundles associated to the chosen metaplectic structure via these modules by $\mathcal{E}^{20}_{\pm},$ $\mathcal{E}^{21}_{\pm}$ and $\mathcal{E}^{22}_{\pm},$ respectively. We define $\mathcal{E}^{2j}:=\mathcal{E}^{2j}_+\oplus \mathcal{E}^{2j}_-$ for  $j=0,1,2.$

In the paper, we shall prove that the symplectic Ricci tensor field maps a symplectic spinor field $\phi \in \Gamma(M, \mathcal{S})$ into $\Gamma(M,\mathcal{E}^{20}\oplus \mathcal{E}^{21})$ and the symplectic Weyl tensor field maps a symplectic spinor field into $\Gamma(M,\mathcal{E}^{21}\oplus \mathcal{E}^{22}).$ For an arbitrary symplectic spinor field $\phi \in \Gamma(M,\mathcal{S}),$ the projections 
of $R^S\phi$ to the invariant spaces $\Gamma(M,\mathcal{E}^{2j})$ ($j=0,1,2$) are explicitly computed.  This describes a structure of the action of the curvature tensor field $R^S$ on the space of symplectic spinor fields in terms of the invariant parts of the underlying connection. In what follows, this result will be called the decomposition result.

The result described above seems to be rather abstract or technical. 
But actually, knowing the decomposition of $R^S\phi$ makes it possible to
derive several conclusions for certain invariant differential operators, which are defined with help of the Fedosov connection.

This is the case of the application we shall mention. Let us briefly describe its context.
In 1994, K. Habermann introduced a symplectic analogue of the
Dirac operator known from Riemannian geometry, the so called symplectic Dirac operator.
The symplectic Dirac operator was introduced with the help of the so
called symplectic Clifford multiplication, see Habermann \cite{Ha}. 
It is possible to define the same operator (up to a complex scalar multiple) using the de Rham sequence tensored (twisted) by symplectic spinor fields as one usually does in the Riemannian spin geometry to get a definition of the Dirac, twistor and Rarita-Schwinger operator and their further  higher spin analogues.
  
  Under the assumption the symplectic Weyl tensor $W$ of the Fedosov connection is trivial, there exists a complex consisting of two differential operators $T_0$ and $T_1.$ These operators will also be called symplectic twistor operators and they will be defined using the de Rham sequence tensored by symplectic spinor fields.  
One of the advantage of the decomposition result is a complete avoidance of possibly lengthy computations in coordinates when proving that $T_0$ and $T_1$ form a complex (provided $W=0$). One can say that the coordinate computations were absorbed into the proof of the decomposition result.
Though finding the complex seems to be a rather particular result,  there is a strong hope of deriving a longer complex under the same assumption.


The reader interested in applications of this theory in physics is
referred to Green, Hull \cite{GH}, where the symplectic spinors are
used in the context of 10 dimensional super string theory. In Reuter
\cite{Reuter}, symplectic spinors are used in the theory of  the so called
Dirac-K\"{a}hler fields.

In the second section, some   basic facts on
 the  symplectic spinor representation and higher symplectic spinors are recalled.
 In the section 3, basic properties of  torsion-free symplectic (Fedosov) connections and their curvature tensor field are mentioned.
 In the fourth section, the action of the  curvature tensor field $R^S$ of the associated  symplectic spinor covariant derivative $\nabla^S$ acting on the space of  symplectic spinor fields described (Corollary 11).
 In this section, the mentioned complex of the two symplectic twistor operators is presented (Theorem 12).

\section{Metaplectic representation, higher symplectic spinors and basic notation}

We start with a summary of notions from representation theory, which we shall need in this paper. From the point of view of this article, the notions are rather of a technical character.  
Let $G$ be a reductive Lie group in the sense of Vogan (see Vogan \cite{Vogan}), $\mathfrak{g}$ be the Lie algebra of $G$ and $K$ be a maximal compact subgroup of $G.$  Typical examples of reductive groups are finite covers of semisimple Lie subgroups of a general linear group of a finite dimensional vector space.
Let $\mathcal{R}(G)$ be the category the object of which are complete, locally convex, Hausdorff topological spaces with continuous linear $G$-action, such that the resulting representation is admissible and of finite length; the morphisms are continuous $G$-equivariant linear maps between the objects. Let $\mathcal{HC}(\mathfrak{g},K)$ be the category of Harish-Chandra $(\mathfrak{g},K)$-modules and let us consider the forgetful Harish-Chandra functor $HC:\mathcal{R}(G)\to \mathcal{HC}(\mathfrak{g},K).$
It is well known that there exists an adjoint functor $mg: \mathcal{HC}(\mathfrak{g},K)\to \mathcal{R}(G)$ to the Harish-Chandra functor $HC$. This functor is usually called the minimal globalization  functor and its existence is a deep result in representation theory. For details and for the existence of the minimal globalization functor $mg,$ see Kashiwara, Schmid \cite{KS} and/or Vogan \cite{Vogan}. 

For a  representation
${\bf E} \in \mathcal{R}(G)$ of $G,$ we shall denote the corresponding Harish-Chandra $(\mathfrak{g},K)$-module $HC({\bf E})$ by $E.$ When we will only  be considering its $\mathfrak{g}$-module structure, we shall use the symbol $\mathbb{E}$ for it. 

Now, suppose $K$ is connected and two complex modules $E,F \in \mathcal{HC}(\mathfrak{g},K)$ are given such that both $\mathbb{E}$ and $\mathbb{F}$ are irreducible highest weight $\mathfrak{g}^{\mathbb{C}}$-modules. Because $mg$ is an adjoint functor to the functor $HC,$ we have $\mbox{Hom}_G(mg(E),mg(F))\simeq \mbox{Hom}_{(\mathfrak{g},K)}(E,F).$ It is well known that the category of $(\mathfrak{g},K)$-modules is a full subcategory of the category of $\mathfrak{g}$-modules provided $K$ is connected. Due to that, we have $\mbox{Hom}_{(\mathfrak{g},K)}(E,F)\simeq \mbox{Hom}_\mathfrak{g}(\mathbb{E},\mathbb{F}).$ Because $\mathbb{E}$ and $\mathbb{F}$ are complex irreducible highest weight modules over $\mathfrak{g}^{\mathbb{C}}$, the Dixmier's version of the Schur lemma implies $\mbox{dim} \mbox{Hom}_{\mathfrak{g}}(\mathbb{E},\mathbb{F})=1$ iff $\mathbb{E} \simeq \mathbb{F}$ (see Dixmier \cite{D}, Theorem 2.6.5 and Theorem 2.6.6). 
Summing up, we have
$\mbox{dim} \mbox{Hom}_{G}(mg(E),mg(F)) = 1$ iff $\mathbb{E}\simeq \mathbb{F}.$  
For brevity, we will refer to this simple statement as to the {\it globalized Schur lemma}.  

Further, if $(p:\mathcal{G} \to M,G)$ is a principal
$G$-bundle, we shall denote the vector bundle associated to this
principal bundle via a representation $\sigma: G \to
\hbox{Aut}(\bf W)$ of $G$ on ${\bf W}$ by $\mathcal{W},$ i.e.,
$\mathcal{W}=\mathcal{G}\times_{\sigma} {\bf W}.$
Let us also mention that we shall often use the Einstein summation convention for repeated indices (lower and upper) without mentioning it explicitly.

Now, we shall focus our attention to the studied case, i.e., to the symplectic one.
To fix a notation, let us recall some notions from the symplectic linear algebra.
Let us consider a real symplectic vector space
$(\mathbb{V},\omega_0)$ of dimension $2l,$ i.e., $\mathbb{V}$ is a
$2l$ dimensional real vector space and $\omega_0$ is a
non-degenerate antisymmetric bilinear form on $\mathbb{V}.$ Let us
choose two Lagrangian subspaces\footnote{maximal isotropic wr. to $\omega_0$} $\mathbb{L}, \mathbb{L}' \subseteq
\mathbb{V}$ such that $\mathbb{L}\oplus \mathbb{L}'=\mathbb{V}.$ It follows that $\mbox{dim}(\mathbb{L})=\mbox{dim}(\mathbb{L}')=l.$
Throughout this article, we shall use a symplectic basis
$\{e_i\}_{i=1}^{2l}$ of $\mathbb{V}$ chosen  in such a way   that
$\{e_i\}_{i=1}^l$ and $\{e_i\}_{i=l+1}^{2l}$ are respective bases of
$\mathbb{L}$ and $\mathbb{L}'.$ Because the definition of  a symplectic
basis is not unique, let us fix   one which   shall be  used in this
text. A basis $\{e_i\}_{i=1}^{2l}$ of $\mathbb{V}$ is called
symplectic basis of $(\mathbb{V},\omega_0)$ if
$\omega_{ij}:=\omega_0(e_i,e_j)$ satisfies $\omega_{ij}=1$ if and
only if $i\leq l$ and $j=i+l;$ $\omega_{ij}=-1$ if and only if $i>l$
and $j=i-l$ and finally, $\omega_{ij}=0$ in other cases. Let
$\{\epsilon^i\}_{i=1}^{2l}$ be the basis of $\mathbb{V}^*$ dual to
the basis $\{e_i\}_{i=1}^{2l}.$  For $i,j=1,\ldots,
2l,$ we define $\omega^{ij}$ by
$\sum_{k=1}^{2l}\omega_{ik}\omega^{jk}=\delta_i^j,$ for $i,j=1,\ldots,
2l.$ Notice that not only $\omega_{ij}=-\omega_{ji},$ but also
$\omega^{ij}=-\omega^{ji},$ $i,j =1, \ldots, 2l.$

Let us denote the symplectic group of $(\mathbb{V},\omega_0)$  by
$G,$ i.e., $G :=Sp(\mathbb{V},\omega_0)\simeq Sp(2l,\mathbb{R}).$
Because the maximal compact subgroup $K$ of $G$ is isomorphic to the
unitary group $K \simeq U(l)$ which is of homotopy type
$\mathbb{Z},$ there exists a nontrivial two-fold covering
$\tilde{G}$ of $G.$ See, e.g., Habermann, Habermann \cite{HH} for details. This two-fold covering is called metaplectic
group  of $(\mathbb{V},\omega_0)$ and it is denoted by
$Mp(\mathbb{V},\omega_0)$. Let us remark that $Mp(\mathbb{V},\omega_0)$ is reductive in the sense of Vogan.
In the considered case, we have
$\tilde{G}\simeq Mp(2l,\mathbb{R}).$ For a later use, let us reserve
the symbol $\lambda$ for the mentioned covering. Thus $\lambda:
\tilde{G} \to G$  is a fixed member  of the isomorphism class of all
nontrivial  $2:1$ coverings of $G$. Because $\lambda:\tilde{G}\to G$
is a homomorphism of Lie groups and $G$ is a subgroup of the general
linear group $GL(\mathbb{V})$ of $\mathbb{V},$ the mapping $\lambda$ is
also a representation of the metaplectic group $\tilde{G}$ on the
vector space $\mathbb{V}.$ Let us define $\tilde{K}:=\lambda^{-1}(K).$ Then $\tilde{K}$ is a maximal compact subgroup of $\tilde{G}.$
 One can easily see that $\tilde{K}\simeq \widetilde{U(l)}:=\{(g,z)\in U(l)\times \mathbb{C}^{\times}| \mbox{det}(g)=z^2\}$ and thus, $\tilde{K}$ is connected. The Lie algebra of $\tilde{G}$ is isomorphic to the Lie algebra $\mathfrak{g}$ of $G$ and we will identify them. One has $\mathfrak{g}=\mathfrak{sp}(\mathbb{V},\omega_0)\simeq \mathfrak{sp}(2l,\mathbb{R}).$  

 From now on, we shall restrict ourselves to the case $l\geq 2$
  without mentioning it explicitly. The case $l=1$ should be handled separately (though analogously) because
  the shape of the root system of $\mathfrak{sp}(2,\mathbb{R})\simeq \mathfrak{sl}(2,\mathbb{R})$ is different from that
  one of
  of the root system of $\mathfrak{sp}(2l,\mathbb{R})$ for $l>1.$
   As usual, we shall denote the
complexification of $\mathfrak{g}$ by $\mathfrak{g}^{\mathbb{C}}.$
Obviously, $\mathfrak{g}^{\mathbb{C}}\simeq
\mathfrak{sp}(2l,\mathbb{C}).$  Let us choose a Cartan subalegbra $\mathfrak{h}^{\mathbb{C}}$ of $\mathfrak{g}^{\mathbb{C}}$ and an ordering on the set of roots of $(\mathfrak{g}^{\mathbb{C}},\mathfrak{h}^{\mathbb{C}}).$
If $\mathbb{E}$ is an irreducible highest weight $\mathfrak{g}^{\mathbb{C}}$-module with a highest weight $\lambda,$ we shall denote it by the  symbol $L(\lambda).$  
Let us denote the fundamental weight basis 
of $\mathfrak{g}^{\mathbb{C}}$ wr. to the above choices by $\{\varpi _{i}\}_{i=1}^l.$


\subsection{Metaplectic representation and symplectic spinors}

There exists a distinguished infinite  dimensional unitary representation of
the metaplectic group $\tilde{G}$ which does not descend to a
representation of the symplectic group $G.$ This representation,
called {\it Segal-Shale-Weil},\footnote{The names oscillator or
metaplectic representation are also used in the literature. We shall
use the name Segal-Shale-Weil  in this text, and reserve the name
metaplectic for certain representation arising from the
Segal-Shale-Weil one.} plays a fundamental role in geometric
quantization of Hamiltonian mechanics, see, e.g., Woodhouse
\cite{Wood}, and in the theory of modular forms and theta
correspondence, see, e.g., Howe \cite{H}. We shall not give a
definition of this representation here and  refer the interested
reader to Weil \cite{Weil} or  Habermann, Habermann \cite{HH}.
 
 The Segal-Shale-Weil representation, which we shall denote by $U$
here, is a complex infinite dimensional unitary representation of
$\tilde{G}$ on the space of complex valued square Lebesgue
integrable functions defined on the Lagrangian subspace
$\mathbb{L},$ i.e.,
$$U: \tilde{G} \to
\mathcal{U}({\bf L^2}(\mathbb{L})),$$ where $\mathcal{U}({\bf W})$
denotes the group of unitary operators on a Hilbert space ${\bf W}.$
In order to be precise, let us refer to the space ${\bf
L^2}(\mathbb{L})$ as to the Segal-Shale-Weil module. It is known  that the Segal-Shale-Weil module belongs to the category $\mathcal{R}(\tilde{G}).$ (See   Kashiwara, Vergne \cite{KV} for details and Segal-Shale-Weil representation in general.)
It is easy to see
that this representation  splits into two
irreducible modules ${\bf L^{2}}(\mathbb{L})\simeq {\bf
L^{2}}(\mathbb{L})_+\oplus {\bf L^{2}}(\mathbb{L})_-.$ The first
module consists of even and the second one of  odd complex valued square
Lebesgue integrable  functions on the Lagrangian subspace
$\mathbb{L}.$ Let us remark that a typical construction of the
Segal-Shale-Weil representation is based on the so called
Schr\"{o}dinger representation of the Heisenberg group of
$(\mathbb{V}=\mathbb{L}\oplus\mathbb{L}',\omega_0)$ and a use of the
Stone-von Neumann theorem.
 
 For technical reasons, we shall need the minimal
globalization of the underlying $(\mathfrak{g},\tilde{K})$-module $HC({\bf L^2}(\mathbb{L}))$  of the introduced Segal-Shale-Weil module.  We
shall call this minimal globalization {\it metaplectic representation} and
denote it by $meta,$ i.e.,
$$meta: \tilde{G} \to \hbox{Aut}(mg(HC({\bf L^2}(\mathbb{L})))),$$ where $mg$ is the minimal globalization functor (see this section and the references therein).
For our convenience, let us denote the module
$mg(HC({\bf L^2}(\mathbb{L})))$ by ${\bf S}.$ Similarly we define $\bf S_+$ and $\bf S_-$
to be the minimal globalizations of the underlying Harish-Chandra modules of the modules
${\bf L^2}(\mathbb{L})_+$ and ${\bf L^{2}}(\mathbb{L})_-$ introduced above.
Accordingly to ${\bf L^{2}}(\mathbb{L})\simeq {\bf L^{2}}(\mathbb{L})_+\oplus
{\bf L^{2}}(\mathbb{L})_-,$  we have $\bf S \simeq {\bf S_+} \oplus {\bf S_-}.$
We shall call the $Mp(\mathbb{V},\omega)$-module   $\bf S$
the  symplectic spinor module  and its elements {\it symplectic spinors}.  For
the name "spinor", see Kostant \cite{Kostant2} or the Introduction.

Further notion related to the symplectic vector space
$(\mathbb{V}=\mathbb{L}\oplus \mathbb{L}',\omega_0)$ is the so called symplectic Clifford
multiplication of elements of ${\bf S}$ by vectors from $\mathbb{V}.$ For a symplectic spinor $f\in {\bf
S},$ we define
$$(e_i.f)(x):=\imath x^i f(x),\footnote{The symbol $\imath$ denotes the imaginary unit, $\imath = \sqrt{-1}.$} $$
$$(e_{i+l}.f)(x):=\frac{\partial f}{\partial x^{i}}(x),
x=\sum_{i=1}^{l}x^i e_i \in \mathbb{L}, i= 1,\ldots,  l.$$ Extending
this multiplication $\mathbb{R}$-linearly, we get the mentioned
symplectic Clifford multiplication.
Let us mention that the multiplication and the differentiation make sense for any $f\in {\bf S}$ because of the interpretation of the minimal globalization. See Vogan \cite{Vogan} for details. Let us remark that in the physical literature, the symplectic Clifford multiplication is usually called the Schr\"odinger quantization prescription.
 
The following lemma is an easy consequence of the definition of the symplectic Clifford multiplication.

{\bf Lemma 1:} For $v,w \in \mathbb{V}$ and $s \in {\bf S},$ we have
$$v.w.s-w.v.s=-\imath \omega_0(v,w)s.$$

{\it Proof.} See Habermann, Habermann \cite{HH}, pp. 11. $\Box$

 Sometimes, we shall write $v.w.s$ instead of $v.(w.s)$ for $v,w \in \mathbb{V}$ and a symplectic spinor $s\in {\bf S}$ and similarly for higher number of multiplying elements. Instead of $e_i.e_j.s,$ we shall write $e_{ij}.s$ simply 
and similarly for expressions with higher number of multiplying elements, e.g., $e_{ijk}.s$ abbreviates $e_i.e_j.e_k.s.$
 
\subsection{Higher symplectic spinors}

In this subsection, we shall present a result on a decomposition
of the tensor product of the metaplectic representation with the first and with the second wedge
power of the representation $\lambda^*: \tilde{G} \to
GL(\mathbb{V}^*)$ of $\tilde{G}$ (dual to the representation
$\lambda$) into irreducible summands.  
Let us reserve the symbol $\rho$ for the
mentioned  tensor product representation of $\tilde{G}$, i.e.,
$$\rho: \tilde{G} \to \hbox{Aut}(\bigwedge ^{\bullet}\mathbb{V}^*\otimes {\bf S})$$
and $$\rho(g)(\alpha\otimes s):=\lambda(g)^{*\wedge r}\alpha\otimes
meta(g)s$$ for $r\in \{0,\ldots, 2l\},$ $g\in \tilde{G},$   $\alpha\in \bigwedge^r\mathbb{V}^*,$  $s\in \bf{S}$ and extend it linearly.  For definiteness, let us equip the tensor product
$\bigwedge^{\bullet}\mathbb{V}^*\otimes {\bf S}$ with the so called
Grothendieck tensor product topology. See Vogan \cite{Vogan} and Treves \cite{Treves} for
details on this topological structure. In a parallel to the Riemannian case, we shall call the elements of $\bigwedge^{\bullet}\mathbb{V}^* \otimes {\bf S}$ {\it higher symplectic spinors}.

In the next theorem,  the modules of the exterior $1$-forms and $2$-forms with values in the module ${\bf S}$ of symplectic spinors are decomposed into irreducible summands. 

{\bf Theorem 2:} For $\mbox{dim}(\mathbb{V}) =:l \geq 2$, the following isomorphisms
$$\mathbb{V}^*\otimes {\bf S} \simeq {\bf E}^{10}_{\pm}\oplus {\bf E}^{11}_{\pm} \, \mbox{ and}$$
$$\bigwedge^2\mathbb{V}^* \otimes {\bf S}_{\pm} \simeq {\bf E}^{20}_{\pm} \oplus  {\bf E}^{21}_{\pm}
\oplus {\bf E}^{22}_{\pm}$$ hold. 
For $j_1=0,1$ and $j_2=0,1,2$ the
modules ${\bf E}^{1j_1}_{\pm}$ and ${\bf E}^{2j_2}$ are uniquely determined by the conditions that first, they are submodules of the corresponding tensor products and second, 
\begin{eqnarray*}
\mathbb{E}^{10}_{+} &\simeq& \mathbb{E}^{20}_{-} \simeq \mathbb{S}_{-} \simeq L(\varpi_{l-1}-\frac{3}{2}\varpi_l), \, \mbox{  }\mathbb{E}^{10}_- \simeq \mathbb{E}^{20}_+ \simeq \mathbb{S}_+ \simeq L(-\frac{1}{2}\varpi_{l}),\\
\mathbb{E}^{11}_+ &\simeq& \mathbb{E}^{21}_- \simeq  L(\omega_1-\frac{1}{2}\omega_l), \, \mbox{   } \mathbb{E}^{11}_-     \simeq \mathbb{E}^{21}_+ \simeq L(\omega_1+\omega_{l-1}-\frac{3}{2}\omega_l),\\
 \mathbb{E}^{22}_+ &\simeq&  L(\omega_{2}-\frac{1}{2}\omega_l) \, \mbox{ and } \, \mathbb{E}^{22}_- \simeq  L(\omega_2+\omega_{l-1}-\frac{3}{2}\omega_l).
\end{eqnarray*}
{\it Proof.} See Kr\'ysl \cite{KryslRarita} or Kr\'ysl \cite{KryslJRT}. $\Box$
 
{\bf Remark:} In this paper, the multiplicity freeness of the previous two decompositions will be used substantially. One can show that the decompositions are multiplicity-free also in the case $l=2.$ (One only needs to modify the prescription for the highest weights of the summands in the decompositions. See Kr\'{y}sl \cite{KryslJRT} for this case). Let us also mention, that the Theorem 2 is a simple consequence of a theorem of  Britten, Hooper, Lemire \cite{BHL}.

Let us set ${\bf E}^{ij_i}:={\bf E}_{+}^{ij_i}\oplus {\bf E}_-^{ij_i},$ for $i=1,2,$ $j_1=0,1$ and $j_2=0,1,2.$ For the mentioned $i, j_i,$ let us consider
the projections $p^{ij_i}:\bigwedge^i \mathbb{V}^*\otimes {\bf S} \to {\bf E}_{i j_i}.$
The definition is correct because of the multiplicity freeness of the decomposition of the appropriate tensor products.
In the paper, we shall need explicit formulas for these projections.In order to find them, let us introduce the following  mappings.

For $r=0,\ldots, 2l$ and $\alpha \otimes s \in \bigwedge^r
\mathbb{V}^*\otimes {\bf S},$  we set
\begin{eqnarray*}
X&:& \bigwedge ^{r}\mathbb{V}^*\otimes {\bf S} \to
\bigwedge ^{r+1}\mathbb{V}^*\otimes {\bf S,} \, \mbox{  }X(\alpha \otimes
s):=-\sum_{i=1}^{2l}\epsilon^i\wedge \alpha \otimes e_i.s \, ;\\
Y&:& \bigwedge ^{r}\mathbb{V}^*\otimes {\bf S} \to \bigwedge ^{r-1}\mathbb{V}^*\otimes {\bf S}, \, \mbox{  }
Y(\alpha \otimes s):=\sum_{i=1}^{2l}\omega^{ij} \iota_{e_i}\alpha \otimes
e_j.s \, \mbox{ and}\\
H&:& \bigwedge ^{r}\mathbb{V}^*\otimes {\bf S} \to \bigwedge ^{r}\mathbb{V}^*\otimes {\bf S}, \,
H :=\{X,Y\}=XY + YX.
\end{eqnarray*}
  
In order to be able to use these operators in a geometric setting, we shall need the following lemma.

{\bf Lemma 3:} The homomorphisms $X,Y,H$ are $\tilde{G}$-equivariant
with respect to the representation $\rho$ of $\tilde{G}.$

{\it Proof.} This can be verified by a direct computation. See Kr\'ysl \cite{KryslRarita} or 
Kr\'ysl \cite{KryslJRT} for a  proof.
$\Box$

In the next lemma, the values of   $H$ on  homogeneous
components of $\bigwedge^{\bullet}\mathbb{V}^*\otimes {\bf S}$ are computed.

{\bf Lemma 4:} Let $(\mathbb{V},\omega_0)$ be a $2l$ dimensional
symplectic vector space. Then for $r=0,\ldots, 2l,$ we have
$$H_{|\bigwedge^r\mathbb{V}^*\otimes {\bf S}}=\imath(r-l)\hbox{Id}_{|\bigwedge^r\mathbb{V}^*\otimes {\bf S}}.$$

{\it Proof.} This can be verified by a direct computation as well. See Kr\'ysl \cite{KryslRarita} or Kr\'ysl \cite{KryslJRT} for a proof. $\Box$
 
In the next lemma, the projections $p^{2j},$ $j=0,1,2,$ are computed explicitly with help of the operators $X$ and $Y.$

{\bf Lemma 5:} For $l>1,$ the following equalities hold on $\bigwedge^2\mathbb{V}^*\otimes {\bf S}.$ 
\begin{eqnarray}
p^{20}&=&\frac{1}{l}X^2Y^2, \\
p^{21}&=&\frac{\imath}{1-l}(XY-\frac{\imath}{l}X^2Y^2) \, \mbox{ and}\\
p^{22}&=&\mbox{Id}_{|\bigwedge^2\mathbb{V}^*\otimes {\bf S}} -\frac{\imath}{1-l}XY+\frac{1}{1-l}X^2Y^2.
\end{eqnarray}

{\it Proof.} 
\begin{itemize}
\item[1.] From the definition of $Y,$ 
the fact that it is $\tilde{G}$-equivariant (Lemma 3) and the Theorem 2, 
we know that $Y^2$ maps $\bigwedge^2\mathbb{V}^*\otimes {\bf S}_{\pm}$ into  ${\bf S}_{\pm}.$ Because $X^2$ is $\tilde{G}$-equivariant (Lemma 3), it maps ${\bf S}_{\pm}$ into a submodule of $\bigwedge^2\mathbb{V}^*\otimes ({\bf S}_+ \oplus {\bf S}_-)$ which is a (possibly empty) direct sum of sumbodules isomorphic to ${\bf S}_{\pm}.$ Regarding the multiplicity-free decomposition structure of $\bigwedge^2\mathbb{V}^*\otimes {\bf S}_{\pm}$ (Theorem 2), we see that $p':=X^2Y^2$ maps $\bigwedge^2\mathbb{V}^*\otimes {\bf S}_{\pm}$ into
${\bf E}^{20}_{\pm}.$
Computing the value of $p'$ on the element $\psi:=\omega_{ij}\epsilon^i\wedge \epsilon^j\otimes s$ for an $s \in {\bf S},$ we find
that $p'\psi=l\psi.$ Using the globalized Schur lemma (see the section 2), we have $p^{20}=\frac{1}{l}X^2Y^2.$

\item[2.] 
As in the $1^{st}$ item, it is easy to see that $p'' := XY(\mbox{Id}_{|\bigwedge^2\mathbb{V}^*\otimes {\bf S}}-\frac{1}{l}X^2Y^2)$ maps $\bigwedge^2\mathbb{V}^*\otimes {\bf S}$ into ${\bf E}^{21}.$ Let us consider a symmetric $2$-vector $\sigma \in \odot^2 \mathbb{V}$ and denote its $(i,j)$-th component wr. to the basis $\{e_i\}_{i=1}^{2l}$ by $\sigma^{ij}.$ Computing the value $p''\psi$ for
$\psi:=\sigma^{ij}\epsilon^k\wedge \epsilon^l \omega_{jl}\otimes e_{ik}.s,$ $s\in {\bf S},$ we get
$p''\psi=\imath (l-1)\psi.$ Using the globalized Schur lemma again, we have  $p^{21}=\frac{\imath}{1-l}p''.$ Using the defining identity $H=XY+YX$ and the Lemma 4, we get the formula for $p^{21}$ written in the statement of the lemma.

\item[3.] The third equation follows from the fact $p^{20}+p^{21}+p^{22}=\mbox{Id}_{|\bigwedge^2\mathbb{V}^*\otimes {\bf S}}$ and the preceding two items.
\end{itemize}
$\Box$

\section{Symplectic curvature tensor field}

After we have finihed the algebraic part of this paper, we shall recall some results of Vaisman in \cite{Vaisman} and of Gelfand, Retakh and Shubin \cite{GSR}.  
Let $(M,\omega)$ be a symplectic manifold and 
$\nabla$ be a symplectic torsion-free affine connection.  Such connections are also called Fedosov connections and were used, e.g., in the so called Fedosov quantization. See Fedosov
\cite{Fedosov} for this use. Let us remark, that there is no uniqueness result for Fedosov connections, which one has in the case of Riemannian manifolds and Riemannian connections.
By symplectic and torsion-free, we  mean $\nabla \omega =0$ and $T(X,Y):=\nabla_XY-\nabla_YX-[X,Y]=0$ for all $X,Y \in \mathfrak{X}(M),$ respectively.
The triple $(M,\omega,\nabla)$ will be called a {\it Fedosov manifold}.
  
 To fix our notation, let us recall the classical definition of the curvature tensor $R^{\nabla}$ of the connection $\nabla,$ we shall be using here. Let 
$$R^{\nabla}(X,Y)Z:=\nabla_X\nabla_Y Z - \nabla_Y \nabla_X Z - \nabla_{[X,Y]}Z$$ for
$X,Y,Z \in \mathfrak{X}(M).$ Let us choose a local symplectic frame $\{e_i\}_{i=1}^{2l}$ on a fixed open subset $U\subseteq M.$
We shall use the following convention. For $i,j,k,l =1,\ldots, 2l,$ we set
\begin{eqnarray}
R_{ijkl}:=\omega(R^{\nabla}(e_k,e_l)e_j,e_i).
\end{eqnarray}
Let us remark that the convention is different from that one used in Habermann, Habermann \cite{HH}. We shall often write expressions in which indices $i,j,k$ or $l$ e.t.c. occur. We will implicitly mean $i,j,k$ or $l$ are running from $1$ to the dimension of the manifold $M$ without mentioning it explicitly.

Obviously, one has
\begin{eqnarray}
R_{ijkl}&=&-R_{ijlk} \mbox{  and} \label{as}\\
R_{ijkl}+R_{iklj}+R_{iljk}&=&0 \, \mbox{ (1$^{st}$ Bianchi identity)} \label{1B}.
\end{eqnarray}
One can also prove the identity
\begin{eqnarray}
R_{ijkl}=R_{jikl} \label{sy}.
\end{eqnarray} See Gelfand, Retakh, Shubin \cite{GSR} for the proof.  

For a symplectic manifold with a Fedosov connection, one has also the following simple consequence of the $1^{st}$ Bianchi identity:
\begin{eqnarray}
 R_{ijkl}+R_{jkli}+R_{klij}+R_{lijk}=0 \, \mbox{ (extended $1^{st}$ Bianchi identity)}. \label{extBR}
\end{eqnarray}
 
From the symplectic curvature tensor field $R^{\nabla}$, we can build the symplectic Ricci curvature tensor field  $\sigma^{\nabla}$ defined by the 
classical formula

$$\sigma^{\nabla}(X,Y):=\mbox{Tr}(V \mapsto R^{\nabla}(V,X)Y)$$ for each $X,Y \in \mathfrak{X}(M)$ (the variable $V$ denotes a vector field on $M$). For the chosen frame and $i,j=1,\ldots, 2l$, we define
$$\sigma_{ij}:=\sigma^{\nabla}(e_i,e_j).$$

Further, let us define
\begin{eqnarray*}
\widetilde{\sigma}^{\nabla}_{ijkl}&:=&\frac{1}{2(l+1)}(\omega_{il}\sigma_{jk}-\omega_{ik}\sigma_{jl}+\omega_{jl}\sigma_{ik}-\omega_{jk}\sigma_{il}+2\sigma_{ij}\omega_{kl}),\\
\widetilde{\sigma}^{\nabla}(X,Y,Z,V)&:=&\widetilde{\sigma}_{ijkl}X^iY^jZ^kV^l\, \mbox{ and}\\
W^{\nabla}&:=&R^{\nabla}-\widetilde{\sigma}^{\nabla}.
\end{eqnarray*}
for  local vector fields $X=X^ie_i,$ $Y=Y^je_j,$ $Z=Z^ke_k$ and $V=V^le_l.$
We will call the tensor field $W^{\nabla}$ the symplectic Weyl curvature tensor field.
These tensor fields were already introduced in Vaisman \cite{Vaisman}. We shall often drop the index $\nabla$ in the previous expressions. Thus we shall often write $W,$ $\sigma$ and $\widetilde{\sigma}$ instead of
$W^{\nabla},$ $\sigma^{\nabla}$ and $\widetilde{\sigma}^{\nabla},$ respectively.

In the next lemma, a symmetry of $\sigma$ and an equivalent definition of $\sigma$ are stated.

{\bf Lemma 6:} The symplectic Ricci curvature tensor field $\sigma$ is symmetric and $$R^{ijkl}\omega_{kl}=2\sigma^{ij}.$$

{\it Proof.} The proof follows from the definition of the symplectic Ricci curvature tensor field and the equation (\ref{sy}). See Vaisman \cite{Vaisman} for a proof. $\Box$

{\bf Remark:} As in the Riemannian geometry, we would like to rise and lower indices.
Because the symplectic form $\omega$ is antisymmetric, we should be more careful in this case. 
For coordinates ${K_{ab\ldots c\ldots d}}^{rs \ldots t \ldots u}$ of a tensor field on the considered symplectic manifold $(M,\omega),$  we denote
the expression $\omega^{ic}{K_{ab\ldots c \ldots d}}^{rs \ldots t}$ by 
${{{K_{ab \ldots}}^{i}}_{\ldots d}}^{rs \ldots t}$ and 
${K_{ab\ldots c}}^{rs \ldots t \ldots u}\omega_{ti}$ by ${{{K_{ab \ldots c}}^{rs\ldots}}_{i}}^{\ldots u}$(similarly for other types of tensor fields).

{\bf Remark:} In Vaisman \cite{Vaisman}, one can find a proof of a statement saying that the space of tensors $R \in \mathbb{V}^{\otimes 4}$ ($\mbox{dim}\mathbb{V}=2l$) satisfying the relations (\ref{as}), (\ref{1B}) and (\ref{sy}) is an $Sp(\mathbb{V},\omega)$-irreducible module if $l=1$ and decomposes into a direct sum of two irreducible $Sp(\mathbb{V},\omega)$-submodules if $l>1.$ 

In the next lemma, two properties of the symplectic Weyl tensor field are described.

{\bf Lemma 7:} The symplectic Weyl curvature tensor field is totally trace-free, i.e., $$W^{ijkl}\omega_{ij}=W^{ijkl}\omega_{ik}=W^{ijkl}\omega_{il}=$$
$$W^{ijkl}\omega_{jk}=W^{ijkl}\omega_{jl}=W^{ijkl}\omega_{kl}=0$$
and the following equation
\begin{eqnarray}
W_{ijkl}+W_{lijk}+W_{klij}+W_{jkli}=0 \, \mbox{  (extended $1^{st}$ Bianchi identity for $W$)} \label{extW}
\end{eqnarray}
holds.

{\it Proof.} The proof is straightforward and can be done just using the definitions of the symplectic Weyl curvature tensor field $W$, the tensor field $\tilde{\sigma}$ and  the Lemma 6. $\Box$

\section{Metaplectic structure and the curvature tensor on symplectic spinors fields}

Let us start describing the geometric structure with help of which the action of the
symplectic curvature tensor field on symplectic spinors, and the symplectic twistor operators  are defined.  This structure, called  metaplectic, is a precise 
symplectic analogue of the notion of a spin structure in  the Riemannian geometry.

For a symplectic manifold $(M^{2l}, \omega)$  of dimension $2l,$
let us denote the bundle of symplectic reperes in $TM$ by
$\mathcal{P}$ and  the foot-point projection of $\mathcal{P}$ onto
$M$ by $p.$ Thus $(p:\mathcal{P}\to M, G),$ where $G\simeq
Sp(2l,\mathbb{R}),$ is a principal $G$-bundle over $M$. As in
the subsection 2, let $\lambda: \tilde{G}\to G$ be a member
of the isomorphism class of the non-trivial two-fold coverings of
the symplectic group $G.$ In particular, $\tilde{G}\simeq
Mp(2l,\mathbb{R}).$ Further, let us consider a principal
$\tilde{G}$-bundle $(q:\mathcal{Q}\to M, \tilde{G})$ over the
symplectic manifold $(M,\omega).$ We call a pair
$(\mathcal{Q},\Lambda)$   metaplectic structure if  $\Lambda:
\mathcal{Q} \to \mathcal{P}$ is a surjective bundle homomorphism
over the identity on $M$ and if the following diagram,
$$\begin{xy}\xymatrix{
\mathcal{Q} \times \tilde{G} \ar[dd]^{\Lambda\times \lambda} \ar[r]&   \mathcal{Q} \ar[dd]^{\Lambda} \ar[dr]^{q} &\\
                                                            & &M\\
\mathcal{P} \times G \ar[r]   & \mathcal{P} \ar[ur]_{p} }\end{xy}$$
with the
horizontal arrows being respective actions of the displayed groups, commutes.
See, e.g.,  Habermann,  Habermann \cite{HH} and Kostant \cite{Kostant2} for
details on   metaplectic structures. Let us only remark, that typical examples of symplectic manifolds admitting a metaplectic structure are cotangent bundles of orientable manifolds (phase spaces), Calabi-Yau manifolds and complex projective spaces $\mathbb{CP}^{2k+1}$, $k \in \mathbb{N}_0.$

Let us denote the  vector bundle
 associated to the introduced principal $\tilde{G}$-bundle
$(q:\mathcal{Q}\to M,\tilde{G})$  via the representation $\rho$ (introduced in the section 2) restricted to ${\bf S}$ by $\mathcal{S}$ and call
 this associated vector bundle {\it symplectic spinor bundle}.  Thus, we have $\mathcal{S}=\mathcal{Q}\times_{\rho}{\bf S}.$ The sections $\phi \in \Gamma(M,\mathcal{S}),$ will be called {\it symplectic spinor fileds}.
 Further for $i=1,2$ and $j_1=0,1$ and $j_2=0,1, 2$, we define the  associated vector bundles $\mathcal{E}^{ij_i}$ by the prescription:  $\mathcal{E}^{ij_i}:=\mathcal{Q}\times_{\rho} {\bf E}^{ij_i}.$   

Because the projections $p^{10}, p^{11}, p^{20},  p^{21}$ and  $p^{22} $ and the operators
$X,Y$ and $H$ are $\tilde{G}$-equivariant (Lemma 3), they lift to operators
acting on sections of the corresponding associated vector bundles.
We shall use the same symbols as for the defined operators as for
their "lifts" to the associated vector bundle structure.

\subsection{Curvature tensor on symplectic spinor fields}

Let $(M,\omega,\nabla)$ be a Fedosov manifold admitting a metaplectic structure $(\mathcal{Q},\Lambda).$ 
The (symplectic) connection $\nabla$ determines the associated principal bundle connection $Z$
on the principal bundle $(p:\mathcal{P}\to M, G).$ This connection lifts to a principal bundle connection on  the principal bundle
$(q:\mathcal{Q}\to M,\tilde{G})$ and defines the associated convariant on the symplectic bundle $\mathcal{S},$ which we shall denote by $\nabla^S$ and call {\it symplectic spinor covariant derivative}. 
The curvature tensor field $R^{S}$ on the symplectic spinor bundle is given by the classical formula
$$R^S:=d^{\nabla^S} \nabla^{S},$$ where $d^{\nabla^S}$ is the associated exterior covariant derivative.
 
 In the next lemma, the action of $R^S$ on the sapce of symplectic spinors is described using just the symplectic curvature tensor field $R.$
 
{\bf Lemma 8:} Let $(M,\omega,\nabla)$ be a 
Fedosov manifold  admitting a metaplectic structure.
Then for a symplectic spinor field $\phi \in \Gamma(M,\mathcal{S}),$ we have
$$R^S\phi=\frac{\imath}{2} {R^{ij}}_{kl}\epsilon^k\wedge \epsilon^l \otimes e_i.e_j.\phi.$$
  
{\it Proof.} See  Habermann, Habermann \cite{HH} pp. 42. $\Box.$
 
Let us define the tensor fields $\sigma^S$ and $W^S$  by the formulas
$$\sigma^S\phi:=\frac{\imath}{2}{\sigma^{ij}}_{kl}\epsilon^k\wedge\epsilon^l\otimes e_i.e_j.\phi \, \mbox{ and}$$ 
$$W^S\phi:= \frac{\imath}{2}{W^{ij}}_{kl}\epsilon^k\wedge \epsilon^l \otimes e_{i}.e_{j}.\phi$$ for a symplectic spinor field $\phi \in \Gamma(M,\mathcal{S}).$

{\bf Theorem 9:} Let $(M,\omega,\nabla)$ be a Fedosov manifold admitting a metaplectic structure.   Then for a symplectic spinor field
$\phi \in \Gamma(M,\mathcal{S}),$ we have $$\sigma^S\phi \in \Gamma(M,\mathcal{E}^{20}\oplus \mathcal{E}^{21}).$$

{\it Proof.} Using the definition  of $\tilde{\sigma}$ and the Lemma 1 repeatedly we have for $\phi \in \Gamma(M,\mathcal{S}),$
\begin{eqnarray*}
\frac{2}{\imath}\sigma^S\phi
&=&\mbox{$\tilde{\sigma}^{ij}$}_{kl}\epsilon^k\wedge \epsilon^l\otimes e_{ij}.\phi\\
&=&({\omega^{i}}_l{\sigma^{j}}_{k}-{\omega^{i}}_{k}{\sigma^{j}}_{l}+{\omega^{j}}_{l}{\sigma^{i}}_{k}
-{\omega^{j}}_{k}{\sigma^{i}}_{l}+2\sigma^{ij}\omega_{kl})\epsilon_k\wedge \epsilon_l\otimes e_{ij}.\phi\\
&=&(-{\sigma^{j}}_{k}\epsilon^k\wedge \epsilon^i + {\sigma^{j}}_{k}\epsilon^i\wedge \epsilon^k-{\sigma^{i}}_{k} \epsilon^k\wedge \epsilon^j
+{\sigma^{i}}_{k}\epsilon^j\wedge \epsilon^k +\\ 
&& + 2\sigma^{ij}\omega_{kl} \epsilon^k\wedge \epsilon^l)\otimes e_{ij}.\phi\\
&=&2{\sigma^{j}}_{k} \epsilon^i\wedge\epsilon^k \otimes e_{ij}.\phi - 2{\sigma^{i}}_{k} \epsilon^k\wedge\epsilon^j \otimes e_{ij}.\phi
+ 2\sigma^{ij}\omega_{kl}\epsilon^k\wedge \epsilon^l\otimes e_{ij}.\phi\\
&=& 2{\sigma^{j}}_{k} \epsilon^i\wedge\epsilon^k \otimes (e_{ij}.+e_{ji}.)\phi+2\sigma^{ij}\omega_{kl}\epsilon^k\wedge \epsilon^l\otimes e_{ij}.\phi\\
&=& 2 \sigma^{jl}\omega_{lk}\epsilon^i\wedge \epsilon^k \otimes (e_{ij}.+e_{ji}.)\phi + 2\sigma^{ij}\omega_{kl}\epsilon^k\wedge \epsilon^l\otimes e_{ij}.\phi\\
&=& 4 \sigma^{jl}\omega_{lk}\epsilon^i\wedge\epsilon^k\otimes e_{ij}.\phi+2\sigma^{ij}\omega_{kl}\epsilon^k\wedge \epsilon^l\otimes e_{ij}.\phi
\end{eqnarray*}

It is straightforward but tedious to verify the next identities:
\begin{eqnarray*}
X^2Y^2(2\sigma^{ij}\omega_{kl}\epsilon^k\wedge\epsilon^l\otimes e_{ij}\phi)&=&2l\sigma^{ij}\omega_{kl}\epsilon^k\wedge\epsilon^l\otimes e_{ij}\phi,\\
X^2Y^2(4\sigma^{jl}\omega_{lk}\epsilon^i\wedge \epsilon^k \otimes e_{ij}.\phi)&=&
2\sigma^{ik}\omega_{jm}\epsilon^j\wedge \epsilon^m \otimes e_{ik}.\phi,\\
XY(2\sigma^{ij}\omega_{kl}\epsilon^k\wedge \epsilon^l\otimes e_{ij}.\phi)&=&
-2\imath \sigma^{ij}\omega_{kl}\epsilon^k\wedge \epsilon^l\otimes e_{ij}.\phi \, \mbox{ and}\\
XY(4\sigma^{jl}\omega_{lk}\epsilon^i\wedge \epsilon^k \otimes e_{ij}.\phi)&=&
4\imath (l-1)\sigma^{jl}\omega_{lk}\epsilon^m\wedge \epsilon^k\otimes e_{mj}.\phi-\\
&&2\imath \sigma^{jl}\omega_{mi}\epsilon^m\wedge \epsilon^i\otimes e_{lj}.\phi.
\end{eqnarray*}

Using the formulas (1) and (2),  we get:
 \begin{eqnarray}
p^{20}\sigma^S\phi &=& \imath \sigma^{ij}\omega_{kl}\epsilon^k\wedge \epsilon^l \otimes (e_{ij}.+\frac{1}{l}e_{ij}.)\phi \, \mbox{ and} \label{ric1}\\
p^{21}\sigma^S\phi &=& \imath\sigma^{ij}\epsilon^k\wedge \epsilon^l (2\omega_{il}\otimes e_{kj}.-\frac{1}{l}\omega_{kl}\otimes e_{ij}.) \label{ric2}\phi.
\end{eqnarray}
Adding these two formulas and comparing them with the result of the computation of $\frac{2}{\imath}\sigma^S \phi$, we get $(p^{20}+ p^{21})\sigma^S \phi=\sigma^S \phi.$ Now, the statement follows.
$\Box$

{\bf Theorem 10:} Let $(M,\omega,\nabla)$ be a Fedosov manifold admitting a metaplectic structure.  Then for a symplectic spinor field $\phi \in \Gamma(M,\mathcal{S}),$ we have $$W^S\phi \in \Gamma(M,\mathcal{E}^{21}\oplus \mathcal{E}^{22}).$$ 

{\it Proof.}  
Let us compute $Y^2W^S\phi$ for a symplectic spinor field $\phi \in \Gamma(M,\mathcal{S}).$
\begin{eqnarray*}
\frac{2}{\imath}Y^2W^S\phi&=&Y(\omega^{nm}{W^{ij}}_{kl}\iota_{e_n}(\epsilon^k \wedge \epsilon^l)\otimes e_{mij}.\phi)\\
&=&Y(\omega^{nm}{W^{ij}}_{kl}(\delta^k_n\epsilon^l-\delta^l_n \epsilon^k) \otimes e_{mij}.\phi)\\
&=&Y((\omega^{nm}{W^{ij}}_{nl}\epsilon^l-{W^{ij}}_{kn}\epsilon^k)\otimes e_{mij}.\phi)\\
&=&2\omega^{nm}Y({W^{ij}}_{nl}\epsilon^l \otimes e_{mij}.\phi)\\
&=&2\omega^{pk}\omega^{nm}{W^{ij}}_{nl}\iota_{e_p}\epsilon_l \otimes e_{kmij}.\phi\\
&=&2\omega^{pk}\omega^{nm}{W^{ij}}_{np}e_{kmij}.\phi = 2W^{ijkl}e_{lkij}.\phi.
\end{eqnarray*}

Now, let us use the extended $1^{st}$ Bianchi identity for the symplectic Weyl curvature tensor field, Eq. (\ref{extW}), i.e.,
\begin{eqnarray*}
W^{ijkl}+W^{jkli} + W^{klij} + W^{lijk} =0.
\end{eqnarray*}

Multiplying this identity by the operator $e_{lkij}.,$ using the relation $e_{ij}.-e_{ji}.=-\imath\omega_{ij}$ (Lemma 1) and the fact that the symplectic Weyl tensor field is totally trace free (Lemma 7), we get the following chain of equations.

\begin{eqnarray*}
&&W^{ijkl}e_{lkij}.+W^{jkli}e_{lkij}.+W^{klij}e_{lkij}.+W^{lijk}e_{lkij}.=0,\\
&&W^{ijkl}e_{lkij}.+W^{jkli}e_{l}(e_{ik}.-\imath \omega_{ki})e_j. + W^{klij}(e_{likj}.-\imath \omega_{ki}e_{lj}.)+\\
&&+W^{lijk}(e_{klij}.-\imath \omega_{lk}e_{ij}.)=0,\\
&&W^{ijkl}e_{lkij}.+W^{jkli}(e_{ilkj}.-\imath \omega_{li}e_{kj}.)+W^{klij}(e_{ilkj}.-\imath \omega_{li}e_{kj}.)+\\
&&+W^{likj}(e_{klji}.-\imath \omega_{ij}e_{kl}.)=0,\\
&&W^{ijkl}e_{lkij}.+W^{jkli}(e_{iljk}.-\imath \omega_{kj}e_{il}.)+W^{klij}(e_{iklj}.-\imath \omega_{lk}e_{ij}.)+\\
&&+W^{lijk}(e_{kjli}. - \imath \omega_{lj}e_{ki}.)=0 \, \mbox{ and}\\
&&3W^{ijkl}e_{lkij}.+W^{klij}(e_{ikjl}.-\imath \omega_{lj}e_{ik}.)=0.
\end{eqnarray*}
Continuing in a similar way, we get $4W^{ijkl}e_{lkij}.=0.$
Summing up, we have $Y^2W^S=0.$
Using the relation $(1)$ for $p^{20},$ we have $p^{20}W^S\phi=0.$ Hence the statement follows. $\Box$

Let us consider a symplectic spinor field $\phi\in \Gamma(M,\mathcal{S}).$
By a straightforward way, we get $XYW^S\phi=2{W^{ijk}}_l\epsilon^m\wedge\epsilon^l\otimes e_{mkij}.\phi.$
Using this result, Theorem 10, definition of $W^S$ and the relations $(2)$ and $(3)$ for $p^{21}$ and $p^{22},$ we get
\begin{eqnarray}
p^{21}W^S\phi&=&\frac{2\imath}{1-l}{W^{ijk}}_l\epsilon^m\wedge\epsilon^l\otimes e_{mkij}.\phi \, \mbox{ and}\label{proj1}\\
p^{22}W^S\phi&=&\frac{\imath}{2}{W^{ij}}_{kl}\epsilon^k\wedge\epsilon^l\otimes e_{ij}.\phi-\frac{2\imath}{1-l}{W^{ijk}}_l\epsilon^m\wedge \epsilon^l \otimes e_{mnij}.\phi \label{proj2}.
\end{eqnarray}

Summing up the preceding two theorems, we have the
 
{\bf Corollary 11:} In the situation described in the formulation of the Theorem 10, we have for a symplectic spinor field $\phi \in \Gamma(M,\mathcal{S})$
\begin{eqnarray*}
p^{20}R^S\phi&=&\imath \sigma^{ij}\omega_{kl}\epsilon^k\wedge \epsilon^l \otimes (e_{ij}.+\frac{1}{l}e_{ij}.)\phi,\\
p^{21}R^S\phi&=&\imath\sigma^{ij}\epsilon^k\wedge \epsilon^l (2\omega_{il}\otimes e_{kj}.
-\frac{1}{l}\omega_{kl}\otimes e_{ij}.)\phi +\\ &&+\frac{2\imath}{1-l}{W^{ijk}}_l\epsilon^m\wedge\epsilon^l\otimes e_{mkij}.\phi \, \mbox{  and} \\
p^{22}R^S\phi&=&\frac{\imath}{2}{W^{ij}}_{kl}\epsilon^k\wedge\epsilon^l\otimes e_{ij}.\phi-\frac{2\imath}{1-l}{W^{ijk}}_l\epsilon^m\wedge \epsilon^l \otimes e_{mkij}.\phi.
\end{eqnarray*}
{\it Proof.} The equations follow from the equations (\ref{ric1}), (\ref{ric2}), (\ref{proj1}) and (\ref{proj2}) and the definitions of $\sigma^S$ and $W^S.$ $\Box$

Now, let us turn our attention to the mentioned application of the decomposition result (Corollary 11).
Let $(M,\omega,\nabla)$ be a Fedosov manifold admitting a metaplectic structure $(\mathcal{Q},\Lambda).$
Then we have the associated  bundles $\mathcal{E}^{ij_i} \to M$ ($i=1,2,$ $j_1=0,1$ and $j_2=0,1,2$) and the symplectic spinor covariant derivative $\nabla^S$ as well as the associated exterior covariant derivative $d^{\nabla^S}$   at our disposal.
Let us introduce the following first order $Mp(2l,\mathbb{R})$-invariant differential operators:
\begin{eqnarray*}
T_0&:& \Gamma(M,\mathcal{S}) \to \Gamma(M,\mathcal{E}^{11}), \, \mbox{ } T_0:=p^{11}\nabla^S \, \mbox{ and}\\
T_1&:& \Gamma(M,\mathcal{E}^{11})\to \Gamma(M,\mathcal{E}^{22}), \, \mbox{ } T_1:=p^{22}d^{\nabla^S}_{|\Gamma(M,\mathcal{E}^{11})}.
\end{eqnarray*}
We shall call these operators {\it symplectic twistor operators}. These definitions are symplectic counterparts of the definitions of twistor operators in  Riemannian spin-geometry.
Using the Corollary 11, we get

{\bf Theorem  12:} Let $(M,\omega,\nabla)$ be a Fedosov manifold admitting a metaplectic structure. Suppose the symplectic Weyl tensor field $W=0.$
Then 
$$\xymatrix{
0\ar[r]& \Gamma(M,\mathcal{S}) \ar[r]^{T_0} & \Gamma(M,\mathcal{E}^{11}) \ar[r]^{T_1} & \Gamma(M,\mathcal{E}^{22})}$$ is a complex of first order differential operators.  
 
{\it Proof.} 
Let us suppose $W = 0.$ Then $p^{22}R^S=0$ (due to the Corollary 11). Using the definition of $R^S,$ we have
$0= p^{22} R^S = p^{22} (d^{\nabla^S}\nabla^S)= p^{22} d^{\nabla^S} (p^{11} + p^{10}) \nabla^S = p^{22}d^{\nabla^S}p^{11}\nabla^S+p^{22}d^{\nabla^S}p^{10}\nabla^{S}.$ According to Kr\'ysl \cite{KryslSVF}, $p^{22}d^{\nabla^S}p^{10}=0.$ Thus we have $0=p^{22}d^{\nabla^S}p^{11}\nabla^S=T_1 T_0,$ giving the statement. 
$\Box$

{\bf Remark:} In Kr\'{y}sl \cite{KryslSVF}, the $Mp(2l,\mathbb{R})$-module
$\bigwedge^{\bullet}\mathbb{V}^*\otimes {\bf S}$ was decomposed into irreducible summands.
Let us denote these irreducible summands by ${\bf E}^{ij}$ (the specification of the indices $i,j$ can be found in the mentioned article or in Kr\'{y}sl \cite{KryslJRT}). Similarly as above, we can introduce the projections $p^{ij}:\bigwedge^i\mathbb{V}^*\otimes {\bf S}\to {\bf E}^{ij}.$ In the mentioned article, we proved that 
$p^{i+1,j}d^{\nabla^S}_{|\Gamma(M,\mathcal{E}^{ik})}=0$ for all appropriate $i,k$ and
$j>k+1$ or $j<k-1.$  In the proof of the preceding theorem, we used this information in the case of $i=1, k=0$ and $j=2.$

\end{document}